\documentclass[a4paper, 12pt]{article}

\usepackage[sort&compress]{natbib}
\bibpunct{(}{)}{;}{a}{}{,} 

\usepackage{amsthm, amsmath, amssymb, mathrsfs, multirow, url}
\usepackage{graphicx} 
\usepackage{ifthen} 
\usepackage{amsfonts}
\usepackage[usenames]{color}
\usepackage{fullpage}
\usepackage{tikz}

\usepackage{xcolor}
\usepackage{enumerate}

\usepackage{soul}

\usepackage{authblk}
\usepackage[normalem]{ulem}
\usepackage{marginnote}
\usepackage{fontawesome5}

\theoremstyle{plain}

\newtheorem{prop}{Proposition}

\newtheorem*{tprin}{Typicality Principle}

\theoremstyle{definition}

\theoremstyle{remark}

\newcommand{\prob}{\mathsf{P}}

\newcommand{\pval}{\mathsf{pval}}

\newcommand{\bin}{{\sf Bin}}
\newcommand{\unif}{{\sf Unif}}
\newcommand{\nm}{{\sf N}}

\newcommand{\gam}{{\sf Gamma}}

\newcommand{\chisq}{{\sf ChiSq}}

\newcommand{\RR}{\mathbb{R}}

\newcommand{\XX}{\mathbb{X}}

\newcommand{\UU}{\mathbb{U}}

\newcommand{\TT}{\mathbb{T}}

\renewcommand{\S}{\mathcal{S}}

\newcommand{\iid}{\overset{\text{\tiny iid}}{\,\sim\,}}
\newcommand{\ind}{\overset{\text{\tiny ind}}{\,\sim\,}}

\newcommand{\prior}{\mathsf{Q}}

\title{The typicality principle and its implications for statistics and data science}

\author[1]{Yiran Jiang}
\author[2]{Zeyu Zhang}
\author[3]{Ryan Martin}
\author[2]{Chuanhai Liu}

\affil[1]{\it \small Department of Biostatistics, Yale University}
\affil[2]{\it \small Department of Statistics, Purdue University}
\affil[3]{\it \small Department of Statistics, NC State University}

\date{\today}

\usepackage{caption}
\usepackage{subcaption}

\begin{document}

\maketitle 

\begin{abstract}
A central focus of data science is the transformation of empirical evidence into knowledge.  As such, the key insights and scientific attitudes of deep thinkers like Fisher, Popper, and Tukey are expected to inspire exciting new advances in machine learning and artificial intelligence in years to come.  Along these lines, the present paper advances a novel {\em typicality principle} which states, roughly, that if the observed data is sufficiently ``atypical'' in a certain sense relative to a posited theory, then that theory is unwarranted.  This emphasis on typicality brings familiar but often overlooked background notions like model-checking to the inferential foreground.  One instantiation of the typicality principle is in the context of parameter estimation, where we propose a new typicality-based regularization strategy that leans heavily on goodness-of-fit testing.  The effectiveness of this new regularization strategy is illustrated in three non-trivial examples where ordinary maximum likelihood estimation fails miserably.  We also demonstrate how the typicality principle fits within a bigger picture of reliable and efficient uncertainty quantification.

\smallskip

\emph{Keywords and phrases:} falsification; goodness-of-fit; inferential model; likelihood; model-checking; prediction. 
\end{abstract}


\section{Introduction}
\label{S:intro}



Data science has captured the attention of researchers and practitioners across all areas of science, business, government, etc.  Like statistics, data science as a discipline is primarily concerned with the transformation of empirical evidence into knowledge about our world, which falls under the umbrella of inductive logic:
\begin{quote}
{\em 
In inductive reasoning we are performing part of the process by which new knowledge is created. The conclusions normally grow more and more accurate as more data are included.  
} \citep[][p.~54]{fisher1935b}
\end{quote}
Despite statistics' 100-year head start, data science has arguably already surpassed it as the leader in applied inductive logic.  This shift from statistics to data science is more than just a change in terminology; it reflects an evolution in our approach to inference, prediction, and decision-making problems, drawing on tools and insights from machine learning, artificial intelligence, and beyond.  As data science continues to advance, it is sure to inspire---and be inspired by---groundbreaking foundational work along lines championed by deep thinkers like Fisher, Karl Popper, and John Tukey.

The complexity inherent in modern data sets means many sources of uncertainty and ambiguity, thus making the analysis and ensuing induction argument highly non-trivial.  As such, insights from epistemology are germane, and the dominant school of thought is falsificationism, as laid out in Karl Popper's {\em Logic of Scientific Discovery}.  His key insight was that, in a sequence of experiments that severely test relevant theories, those theories that withstand this scrutiny will have ``proved their mettle'' \citep[][p.~10]{popper1959}, a necessary, but not sufficient, condition for any theory to be called {\em not-false}.  It is only in a limiting sense, as the number/severity of the tests increases, that a theory can earn the not-false distinction.  A challenge in modern empirical sciences, not found in the ``all swans are white'' style examples often considered in philosophy texts, is that the empirical data cannot logically contradict any legitimate theory, so there will inevitably be uncertainty in drawing inferences.  This necessitates the reliable quantification of said uncertainty and this is the focus of our present contribution.  

Behind the falsificationist perspective is the understanding that experiments tend to produce data that are {\em typical}, i.e., look like what is expected in the real world.  Consequently, if the observed data are atypical---or look sufficiently different from what is expected---relative to a posited theory, then it is fair to consider that theory falsified.  This is our proposed {\em typicality principle} in its basic form.  But what does it mean for data to ``look like'' what is expected?\footnote{Strictly speaking, in practical applications, there's a sense in which all data are ``atypical.'' For instance, in models that involve absolutely continuous distributions, all possible data realizations have zero probability under such models.  Fortunately, strict (a)typicality is not necessary for reliable inference---meaningful notions of (a)typicality can be developed with a little creativity.}  It is common to interpret the model's likelihood function as a measure of its quality of fit to the observed data and, in turn, it is common to judge whether data ``look like'' what is expected by the magnitude of the likelihood.  The {\em law of likelihood} \citep[e.g.,][]{edwards1992, hacking.logic.book} formalizes this.  In some cases, however, the likelihood function could be large as a result of some degeneracy, not because the data ``look like'' what is expected under the corresponding theory.  This highlights a shortcoming in the likelihood-centric approach to inductive inference and a need for new perspectives.  The familiar regularization strategies depend only on the posited theory---e.g., penalize theories incompatible with assumed ``sparsity''---and not on the data itself, hence cannot repair the aforementioned flaw on their own.  Alternatively, our notion of typicality focuses on {\em fit} in a nonparametric goodness-of-fit sense as opposed to a parametric model-based, large-likelihood sense.  The typicality principle advanced here is inspired by Tukey's deep insights on model building/checking \citep{tukey.eda.book, tukey1962future}.  
While philosophical principles tend to be ``top down,'' i.e., instructions passed down from a higher authority, Tukey's brand of philosophy is ``bottom up'' \citep[e.g.,][]{dempster2002.tukey, tukey.part3} and, hence, Popperian in spirit.  Indeed, in model-building, for example, no candidate model is God-given, so a model's merits must be earned by providing a satisfactory explanation of the variability in the observed data.  Our proposed typicality principle aims to build on these key ideas for regularized estimation and beyond.


Building on the typicality principle, our main methodological innovation here is a new brand of {\em typicality-focused} regularization strategies.  More specifically, we recommend adoption of the familiar penalized likelihood framework but with a twist: rather than penalizing theories that are incompatible with prior knowledge, we penalize those theories for which the data could be judged atypical, thus ensuring our derived procedures reward theories that fit the data well, aiding interpretation and enhancing efficiency.  We evaluate the performance of our proposed typicality-focused regularization by applying it to several challenging problems that have historically served as points of contention in the foundations of statistics. Our results demonstrate its efficiency in point estimation and uncertainty quantification more broadly, highlighting its potential to address some of the deepest unresolved issues in statistical science. Beyond its practical applications, we delve into the theoretical underpinnings of the typicality principle, uncovering connections---or lack thereof---to other familiar statistical principles. These connections underscore the broader importance of typicality, not only as a methodological tool but also as a conceptual bridge linking various aspects of statistical reasoning.  By positioning the typicality principle within this rich theoretical and applied context, this paper lays the groundwork for future exploration of its implications in data science and beyond.

The remainder of the paper is organized as follows.  Section~\ref{S:setup} sets the context of our discussion and introduces some key concepts and notation.  Section~\ref{S:first} introduces a first basic version of the typicality principle, and our discussion there focuses on statistical intuitions and philosophical considerations.  In the context of parameter estimation, a particular instantiation of the typicality principle results in our new typicality-based regularization strategy which is also detailed there.  The performance of our proposed typicality-based regularization strategy is investigated in Section~\ref{S:applications} in the context of three nontrivial, paradox-laden examples: a mixture model due to Le~Cam, the Neyman--Scott problem, and Stein's mean vector length.  The challenges faced in these three examples are also common in modern data science applications where there is a risk of over-fitting, so our contributions here extend beyond the simple parametric models considered here.  Section~\ref{S:broader} digs deeper, advancing a formal typicality principle and showing how this fits into a general framework that offers provably reliable uncertainty quantification, beyond point estimation, hypothesis testing, etc.  Connections to other statistical principles are offered and a numerical illustration demonstrates the validity and efficiency of the proposed framework in a challenging marginal inference problem.  We conclude in Section~\ref{S:discuss} concludes with a few remarks.

\section{Problem setup}
\label{S:setup}

To set the scene, we introduce our notation, model assumptions, and objectives.  Let $X$, taking values in the sample space $\XX$, denote the observable data, which can be a scalar, a vector, a matrix, or something else.  We posit a statistical model $\{ \prob_\theta: \theta \in \TT \}$ for $X$, which consists of a collection of probability distributions supported on $\XX$, indexed by a parameter $\theta$ in the parameter space $\TT$.  This model could be the kind of low-dimensional models presented in textbooks, e.g., $\prob_\theta = \bin(n,\theta)$ or $\prob_\theta = \gam(\alpha,\beta)$ with $\theta=(\alpha,\beta)$, or the model could be high-dimensional as in the case of a neural network, so the parameter $\theta$ can be scalar-, vector-, or even function-valued.  The distributions $\prob_\theta$ included in the model have associated probability density or mass functions, denoted by $p_\theta$, and for each $x \in \XX$, the function $L_x(\theta) = p_\theta(x)$ is the likelihood function.  

We will assume that there exists a ``true value'' of the uncertain model parameter $\theta$, which we will denote by $\Theta$; cases in which the true distribution does not belong to the posited model can also be considered, but we will not consider such cases here; see \citet{jiang2022estimation} for further discussion along these lines.  Of course, data $X=x$ is observed, and our objective is to make inference on the uncertain $\Theta$, relative to the posited statistical model.  The most basic inferential objective is point estimation, which amounts to selecting a most likely or plausible value of the parameter based on observation $x$; this will be discussed in more detail in Section~\ref{S:first}.  Beyond point estimation, a goal is to quantify uncertainty about $\Theta$.  More specifically, uncertainty associated with relevant hypotheses ``$\Theta \in H$'' is quantified with a data-dependent, numerical degree of plausibility and, if a formal test was desired, hypotheses determined to be sufficiently implausible would be rejected on those grounds.  
Henceforth, for simplicity of presentation, both $H$ and ``$\Theta \in H$'' will refer to hypotheses about the uncertain $\Theta$.  Given the plausibility assignments to hypotheses $H \subseteq \TT$, these can be inverted to find, e.g., relatively small hypotheses $H$ whose complements are relatively implausible.  This brings confidence set construction under our general inference umbrella. 


We will further assume that prior information about $\Theta$ is vacuous, so the aforementioned degrees of belief/confidence cannot be obtained via application of Bayes's theorem with a genuine or informative prior distribution for $\Theta$.  Consequently, achieving our top-level goal of quantifying uncertainty about $\Theta$ requires considerations beyond the familiar Bayesian and fiducial-like frameworks; details are presented in Section~\ref{S:broader} below.

\section{Typicality principle}
\label{S:first}

Here we adopt Popper's falsificationist view, that is, a hypothesis $H$ concerning the uncertain $\Theta$ cannot be directly verified/confirmed based on data $x$, it can only be refuted or not depending on whether or not truthfulness of $H$ and data $x$ are sufficiently contradictory.  Importantly, falsificationism is really our only option: employing, e.g., the opposing Carnapian--Jeffreysian--Jaynesian confirmationist view \citep[e.g.,][]{carnap1962, jeffreys1961, jaynes2003} in a scientifically and mathematically rigorous way requires a genuine prior probability distribution and application of Bayes's theorem, which is out of reach because we assume vacuous prior information.\footnote{``Vacuous prior information'' technically means that we lack the information required to eliminate any candidate prior distributions for $\Theta$ \citep{martin.partial}, hence, all priors and all corresponding Bayesian solutions are plausible.  Since it is impossible to learn about $\Theta$, given $x$, from a Bayesian analysis under such an extreme credal state, the confirmationist solution is strictly out of reach.  But there is a middle-ground between the extreme ``all priors'' and ``a single prior'' credal states and, in such cases, the {\em generalized Bayes} solution advanced in \citet{walley1991} is available.}  As discussed in Section~\ref{S:intro}, our falsificationist perspective suggests the specification of a strategy for assessing if---and, if so, in what sense and to what extent---the data $x$ is ``typical'' relative to a stated hypothesis $H$ about $\Theta$.  Given such an assessment, inference is at least conceptually straightforward.   

We start here with an informal and intuitive statement of the key principle, one that focuses on simple, singleton hypotheses.  Even this intuitive version has important implications when it comes to point estimation.  More formality is given in Section~\ref{S:broader}. 

\begin{tprin}[intuitive version]
If data $x$ is sufficiently atypical relative to the posited model with parameter $\theta$, then the hypothesis $H=\{\theta\}$ is unwarranted. 
\end{tprin}

The thoughtful reader is surely tempted to make a comparison between the typicality principle as stated above and more familiar notions involving the likelihood.  For instance, the law of likelihood \citep[e.g.,][p.~30]{edwards1992} states, roughly, that data $x$ supports the hypothesis $H=\{\theta\}$ less than hypothesis $H'=\{\theta'\}$ if $L_x(\theta) < L_x(\theta')$ and, consequently, if the gap was sufficiently wide, then hypothesis $H$ would be unwarranted.  This in turn leads naturally to the principle of maximum likelihood \citep[e.g.,][]{stigler2007, aldrich1997, fisher1922, fisher1925}, which suggests using the maximally supported singleton hypothesis, i.e., the value $\hat\theta(x) = \arg\max_\theta L_x(\theta)$ with maximal likelihood, as a point estimator for $\Theta$ given $x$.  It's well-known, however, that the likelihood function can be large at $(x,\theta)$ values even though $x$ is highly atypical of a sample from $\prob_\theta$.  This is often the result of some inadequacy of the model, such as degeneracy or overfitting, and commonly manifests as inconsistency: as the information in the sample increases, the maximum likelihood estimator fails to converge in probability to $\Theta$.  Despite the ubiquity of maximum likelihood estimation in textbooks and in applications, and their asymptotic efficiency in regular cases, their general inadequacy has been widely documented.  Indeed, in 1960, Lucian Le~Cam wrote \citep[][Sec.~11]{vaart.lecam.2002}
\begin{quote}
{\em The author is firmly convinced that a recourse to maximum likelihood is justifiable only when one is dealing with families of distributions that are extremely regular.  The cases in which [maximum likelihood] estimates are readily obtainable and have been proved to have good properties are extremely restricted.}
\end{quote}
and, later, in his book \citep{lecam.1986.book}, 
\begin{quote}
{\em The terms ``likelihood'' and ``maximum likelihood'' seem to have been introduced by RA Fisher who seems also to be responsible for a great deal of the propaganda on the merits of the maximum likelihood method... In view of Fisher's vast influence, it is perhaps not surprising that the presumed superiority of the method is still for many an article of faith promoted with religious fervor.  This state of affairs remains, in spite of a long accumulation of evidence to the effect that maximum likelihood estimates are often useless or grossly misleading.}  
\end{quote}
Examples highlighting the inadequacies of maximum likelihood estimators are presented in Section~\ref{S:applications} below. Le~Cam's point here is just that maximum likelihood works fine in some cases but not in others.  That an estimator works well in some cases and not in others on its own creates no foundational concerns.  If, however, a core {\em principle}---the principle of maximum likelihood---itself is untrustworthy, then that is a serious foundational concern: data science isn't a science if we don't have trustworthy principles.  The typicality principle aims to fill this trustworthiness gap.

One interpretation of the maximum likelihood estimator $\hat\theta(x)$ is that it represents the member of the posited statistical model that is closest (in the sense of Kullback--Leiber divergence (\citeauthor{lehmann1983theory}, \citeyear{lehmann1983theory}; \citeauthor{Pardo2018}, \citeyear{Pardo2018})) to the empirical distribution of the data.  Since the latter empirical distribution necessarily represents the data, it might come as a surprise to hear that $x$ might not look like a typical sample from the ``best approximation'' $\prob_{\hat\theta(x)}$ in the posited model.  
The point is that the maximum likelihood estimator depends only on some summary features of the full data, rarely on the full data itself; on the other hand, model building/checking processes use all aspects of the data.  As an intentionally oversimplified example, suppose that $x$ consists of what is assumed to be an iid sample of size $n$ from a normal distribution with unit variance.  Then the maximum likelihood estimator is the average of those values in $x$ and, hence, it only depends on this summary.  So, of course, if two data sets $x$ and $y$ are of the same size and have the same sample means, then the likelihood functions and hence the maximum likelihood estimators are identical.  But if the sample mean of $x$ is 0 and all of the entries in $y$ are 0, then $x$ might look like a typical sample from $\nm(0,1)$ but $y$ definitely does not.  The catch is that the law of likelihood, etc.~all include the caveat {\em relative to the given model} and, since it would be absurd to model a data set $y$ whose entries are all identically 0 by a normal distribution, this potential issue is rarely discussed.  In fact, the above point boils down to the sufficiency principle and, in the right context, is considered by most (including us) to be a desirable feature of likelihood-based inference.  In real applications, however, where data and model are more complex than in the toy example above, it may not be clear that $\hat\theta(x)$ only depends on an unacceptably incomplete summary of $x$ and, hence, that $x$ might not resemble a typical sample from $\prob_{\hat\theta(x)}$; see Section~\ref{S:applications}.  Moreover, the trope ``all models are wrong but some are useful'' suggests that we should anticipate some degree of model misspecification and manually build in some direct model-fit considerations rather than fully trust the data reductions inherent in likelihood-based methods.  

Our proposed implementation of the typicality principle is via regularization.  This will be anchored around the likelihood function, since this is efficient in regular cases, but our proposed regularization differs from that common in the data science literature in some important ways.  Specifically, consider the objective function 
\begin{equation}
\label{eq:obj}
\rho_\lambda^\text{typ}(x,\theta) = \ell_x(\theta) - \lambda \, r_x^\text{typ}(\theta), \quad \theta \in \TT,
\end{equation}
where $\ell_x = \log L_x$ is the usual log-likelihood function, $\lambda \geq 0$ is a scalar tuning parameter, and $r_x^\text{typ}$ is a data-dependent, typicality-encouraging penalty function.  The goal, of course, is to maximize this objective function at a fixed $x$ to obtain a point estimator $\check\theta(x)$, which implicitly depends on the user's choice of $\lambda$ and of the penalty function $r_x^\text{typ}$. 

A noteworthy point is that the penalty $r_x^\text{typ}$ depends on data $x$.  Almost exclusively,\footnote{One exception is the use of {\em empirical} or {\em data-driven} prior distributions proposed and developed in \citet{martin.mess.walker.eb}, \citet{martin.walker.deb}, and elsewhere.} the penalty functions used in the literature do not depend on data---their goal is to push the maximum likelihood estimator towards where $\Theta$ is believed to be based on {\em a priori} knowledge about the problem.  Note that these common penalty functions do not address the question of whether $x$ is a ``typical'' sample from $\prob_\theta$, and nor can they; any measure that is designed to quantify typicality must depend on both $x$ and $\theta$.  

Since the focus is on determining if $x$ is atypical relative to $\prob_\theta$, Tukey's insights on model-checking suggest considering some variation on a goodness-of-fit test in the construction of $r_x^\text{typ}$.  Indeed, throughout this paper, we consider 
\begin{equation}
\label{eq:reg}
r_x^\text{typ}(\theta) = -\log \pval_x(\theta), \quad \theta \in \TT, 
\end{equation}
where $\pval_x(\theta)$ is the p-value associated with a test of how well $\prob_\theta$ fits the data $x$.  A fairly general formulation is as follows.  Suppose $x=(x_1,\ldots,x_n)$ consists of iid samples where $\prob_\theta$ determines a continuous distribution function $F_\theta$ for the components.  If $x$ is a ``typical'' sample from $\prob_\theta$, then $F_\theta(x_1),\ldots,F_\theta(x_n)$ should resemble an iid $\unif(0,1)$ sample, and we can take $\pval_x(\theta)$ to be the p-value associated with, say, the Kolmogorov--Smirnov test of this uniformity.  See \citet{liu2023reweighted} and \citet{jiang2022estimation} for more on the Kolmogorov--Smirnov p-value for assessing what we refer to here as typicality.  In certain special cases, however, other simpler goodness-of-fit assessments may be available.  For example, in Gaussian models where ``residual sums of squares'' are expected to have a suitable chi-square distribution, the p-value associated with such a chi-square test can be used to construct the penalty function $r_x^\text{typ}$; see Section~\ref{SS:neyman.scott}.  The examples in Section~\ref{S:applications} highlight how the penalty \eqref{eq:reg} leads to desirable regularization, correcting for systematic biases present in the maximum likelihood estimators in sufficiently non-regular models.

\section{Examples: non-regular estimation}
\label{S:applications}

In this section we consider three noteworthy examples involving ``non-regular'' models.  The key characteristic shared by these models is that the likelihood function associated with the relevant feature of $\Theta$ is off-target to the extent that the maximum likelihood estimator is nonsense or at least inconsistent.  While the examples below are relatively simple, the non-regularity means that they have certain aspects in common with modern applications involving complex, high-dimensional models.  Here we apply the general typicality-encouraging regularization strategy in \eqref{eq:obj} above.

\subsection{Le~Cam's mixture}

\newcommand{\p}{\alpha}

Consider the normal mixture model introduced in \cite{le1990maximum}. Let $\p$ be a very small known number, such as $10^{-50}$, and consider the model 
\[ \prob_\theta = (1-\p) \, \nm(\mu, 1) + \p \, \nm(\mu, \sigma^2), \]
indexed by $\theta = (\mu, \sigma^2) \in \RR \times (0,\infty)$.  The goal is to estimate the uncertain true value $\Theta = (M,\Sigma^2)$ of the model parameter.  The challenge, of course, is that since $\p$ is so small, the data are not expected to be particularly informative about the value of $\Sigma^2$.  


Given $n$ independent samples $x=(x_1, ..., x_n)$, the likelihood function is 
\[ L_x(\mu, \sigma^2)
 \propto \prod_{i=1}^n \Bigl[ (1-\p) \exp\Bigl\{-\frac{(x_i-\mu)^2}{2} \Bigr\} + \p \sigma^{-1} \exp\Bigl\{-\frac{(x_i-\mu)^2}{2\sigma^2} \Bigr\} \Bigr]. \]
The difficulty with maximum likelihood estimation in this application is that the likelihood function is unbounded; this implies that the maximum likelihood estimator is on the boundary of the parameter space and, hence, does not exist in a practical sense.  To see this, set $\mu$ equal to any of the observed data points, say, $x_1$, and consider the likelihood as a function of $\sigma^2$ alone:
\begin{align*}
L_x(x_1, \sigma^2) 
& \propto \prod_{i=1}^n \Bigl[ (1-\p) \exp\Bigl\{-\frac{(x_i-x_1)^2}{2} \Bigr\} + \p \sigma^{-1} \exp\Bigl\{-\frac{(x_i-x_1)^2}{2\sigma^2} \Bigr\} \Bigr] \\
& = \{ (1-\p) + \p \sigma^{-1} \} \prod_{i=2}^n \Bigl[ (1-\p) \exp\Bigl\{-\frac{(x_i-x_1)^2}{2} \Bigr\} + \p \sigma^{-1} \exp\Bigl\{-\frac{(x_i-x_1)^2}{2\sigma^2} \Bigr\} \Bigr].
\end{align*}
The second factor, i.e., the product over $i=2,\ldots,n$, is bounded as a function of $\sigma^2$; the first factor, $(1-\p) + \p \sigma^{-1}$, is unbounded as $\sigma^2 \to 0$.  This proves the non-existence of the maximum likelihood estimator (in the interior of the parameter space) and highlights a clear shortcoming in the principle of maximum likelihood.  

\begin{figure}[p]
     \centering
     \begin{subfigure}[b]{0.75\textwidth}
         \centering
         \includegraphics[width=\textwidth]{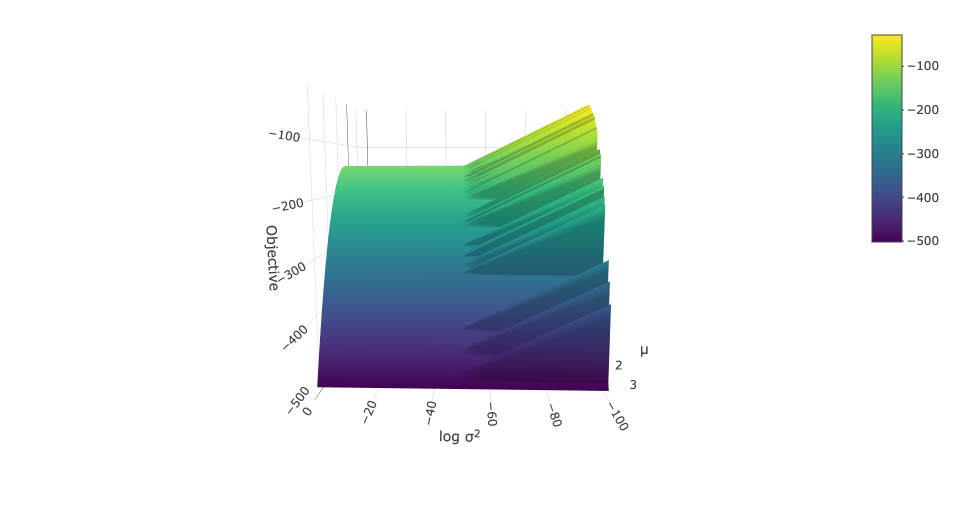}
         \caption{Effect of $\sigma^2$}
         \label{fig:lik1}
     \end{subfigure}
     \\ 
     \begin{subfigure}[b]{0.75\textwidth}
         \centering
         \includegraphics[width=\textwidth]{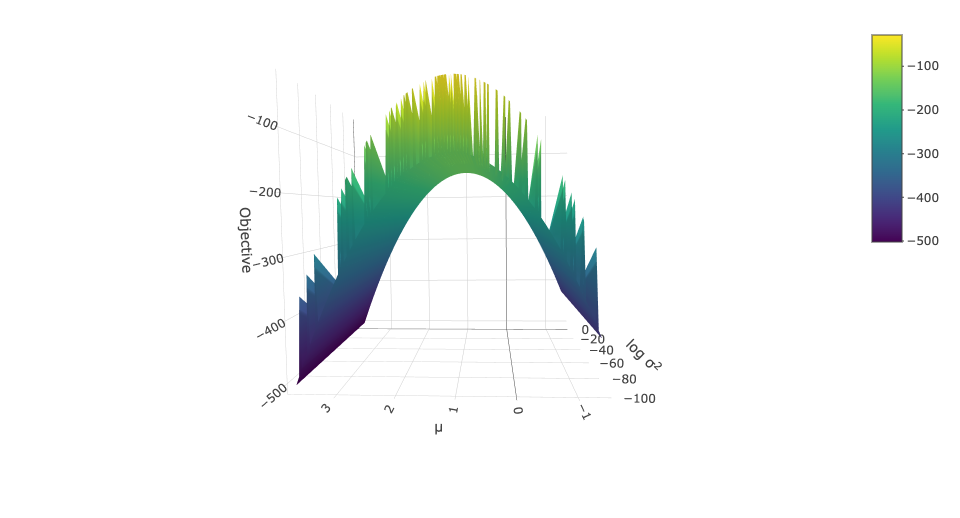}
         \caption{Small $\sigma^2$}
         \label{fig:lik2}
     \end{subfigure}
     \\ 
     \begin{subfigure}[b]{0.75\textwidth}
         \centering
         \includegraphics[width=\textwidth]{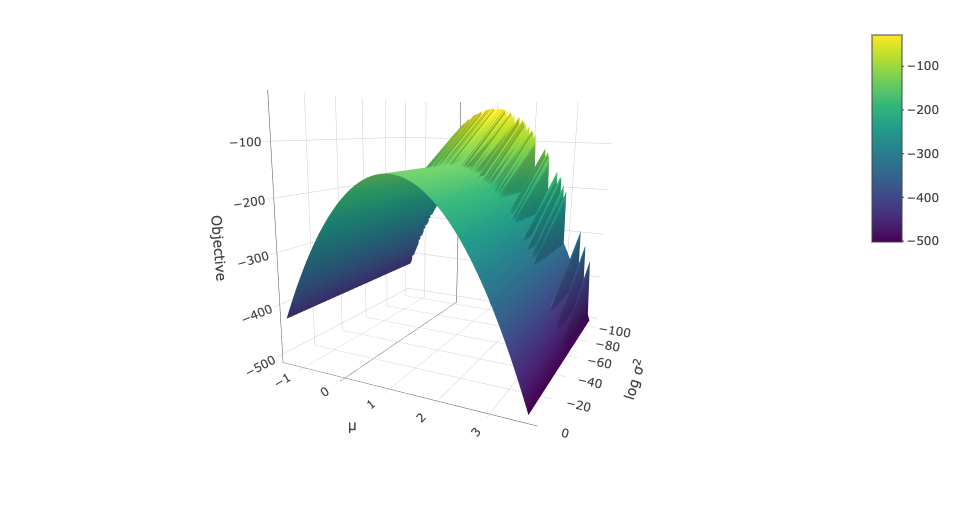}
         \caption{Large $\sigma^2$}
         \label{fig:lik3}
     \end{subfigure}
        \caption{Log-likelihood surfaces from three viewing angles.}
        \label{fig:lik}
\end{figure}

To illustrate this shortcoming, we simulate data of size $n=100$ from the distribution $\prob_\Theta$, where $\Theta=(M,\Sigma^2) = (1,2)$. Then the log-likelihood is evaluated over a range of different $(\mu ,\log\sigma^2)$ values and the results are plotted in Figure~\ref{fig:lik}. Figure~\ref{fig:lik}(\subref{fig:lik1}) shows that $\sigma^2$ will start to dominate $\p$ when it is small enough, driving the log-likelihood to increase monotonically. When $\sigma$ is small enough, as shown in Figure~\ref{fig:lik}(\subref{fig:lik2}), the choices of $\mu$ affect the log-likelihood and create the jumps when they are equal to the values of the samples. On the other hand, the surface is quadratic when $\sigma^2$ is large, as shown in Figure~\ref{fig:lik}(\subref{fig:lik3}), while only the first quadratic plays a role with small $\p$. In summary, using a ``maximum likelihood estimator'' $\hat{\theta}(x) = (x_1, 0)$ amounts to ignoring virtually all of the relevant information in the data, which is very bad practice.  When interest is in making predictions based on the parameter $\hat\theta(x)$, the estimated model $\prob_{\hat\theta(x)}$ will be in efficient, since it is effectively the same fitted model as that based on the single observation $x_1$.

\begin{figure}[p]
     \centering
     \begin{subfigure}[b]{0.75\textwidth}
         \centering
         \includegraphics[width=\textwidth]{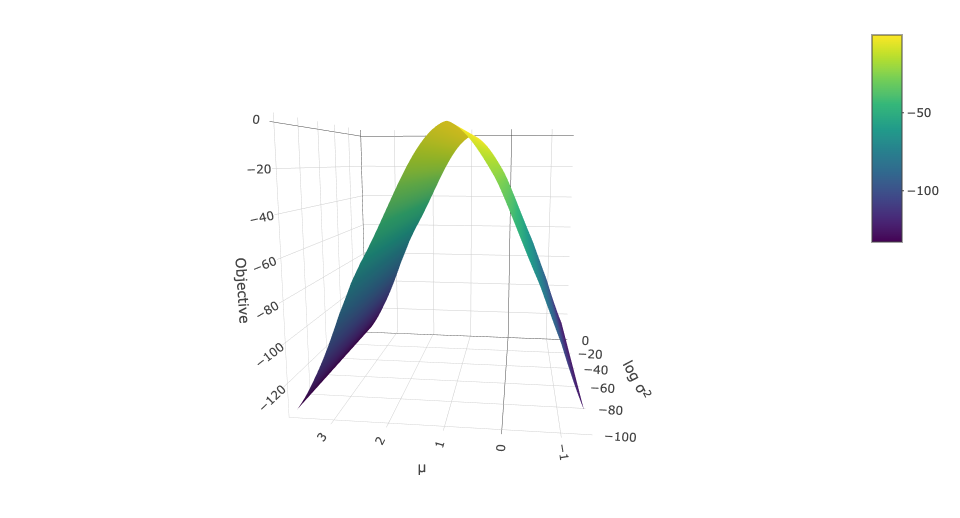}
         \caption{Effect of $\mu$ }
         \label{fig:gg1}
     \end{subfigure}
     \\ 
     \begin{subfigure}[b]{0.75\textwidth}
         \centering
         \includegraphics[width=\textwidth]{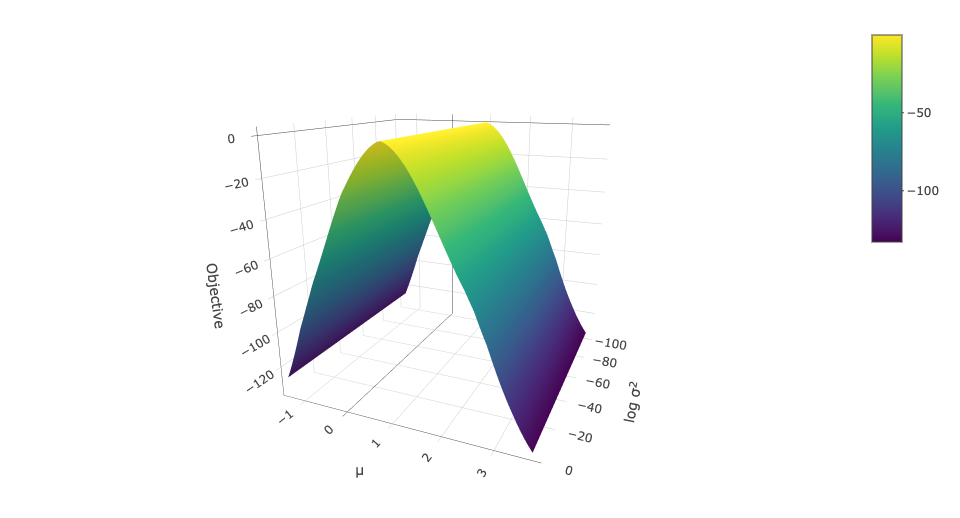}
         \caption{Effects of $\mu$ and $\sigma^2$}
         \label{fig:gg2}
     \end{subfigure}
     \\ 
     \begin{subfigure}[b]{0.75\textwidth}
         \centering
         \includegraphics[width=\textwidth]{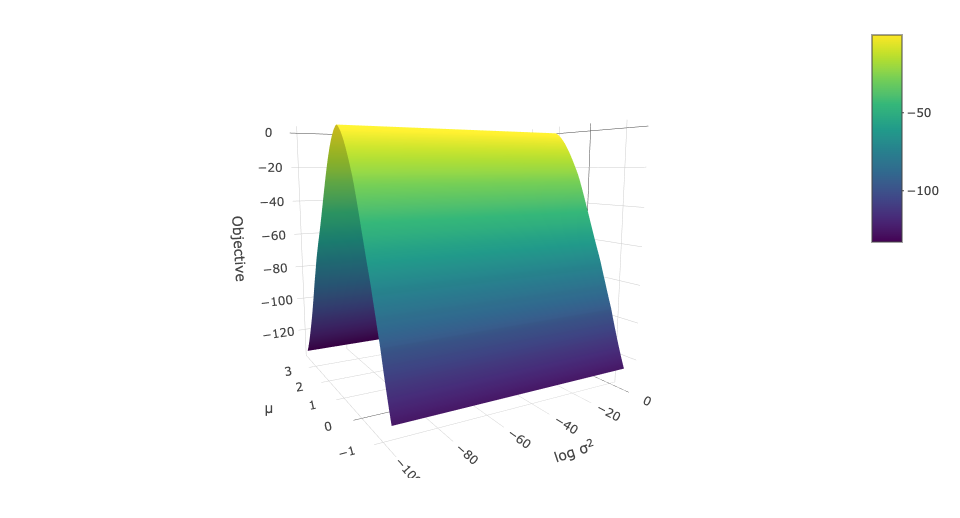}
         \caption{Effect of $\sigma^2$}
         \label{fig:gg3}
     \end{subfigure}
        \caption{Log p-value surfaces for the Kolmogorov--Smirnov test.}
        \label{fig:gg}
\end{figure}

\begin{figure}[t]
     \centering
     \begin{subfigure}[b]{0.475\textwidth}
         \centering
         \includegraphics[width=\textwidth]{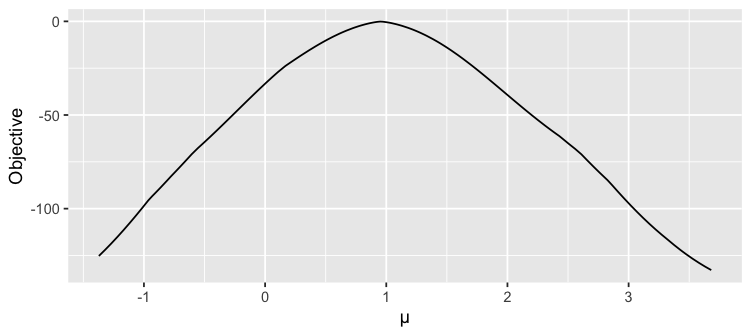}
         \caption{Effect of $\mu$}
         \label{fig:2d1}
     \end{subfigure}
     \hfill
     \begin{subfigure}[b]{0.475\textwidth}
         \centering
         \includegraphics[width=\textwidth]{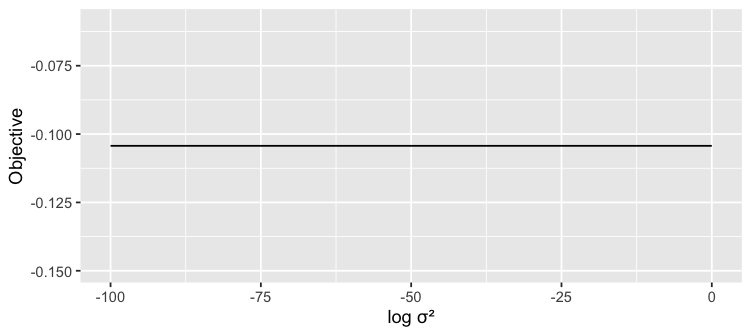}
         \caption{Effect of $\sigma^2$}
         \label{fig:2d2}
     \end{subfigure}
        \caption{Profiled log p-values for the Kolmogorov--Smirnov test.}
        \label{fig:2d}
\end{figure}

To overcome the pitfalls of maximum likelihood, we apply the typicality principle as suggested in Section~\ref{S:first} using the proposed typicality-based regularization \eqref{eq:reg}.  That is, we take $r_x^\text{typ}(\theta)$ as the negative log p-value of the Kolmogorov--Smirnov test assessing the fit of $\prob_\theta$ to the observed $x$.  Clearly from Figure~\ref{fig:gg}, we find that $\mu$ plays a more important role than $\sigma^2$. The choice of $\sigma^2$, however, does not affect the p-values at all.  The p-value profiles, i.e., fixing one parameter and maximizing over the other, are shown in Figure~\ref{fig:2d}. Interestingly, the value of $\mu$ suggested by the goodness-of-fit test is around 1, which coincides with the true parameter value, $M=1$.


\begin{figure}[t]
     \centering
     \begin{subfigure}[b]{0.475\textwidth}
         \centering
         \includegraphics[width=\textwidth]{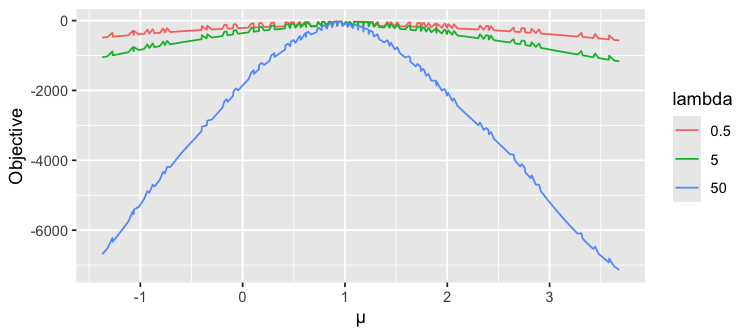}
         \caption{Effect of $\mu$}
         \label{fig:trust1}
     \end{subfigure}
     \hfill
     \begin{subfigure}[b]{0.475\textwidth}
         \centering
         \includegraphics[width=\textwidth]{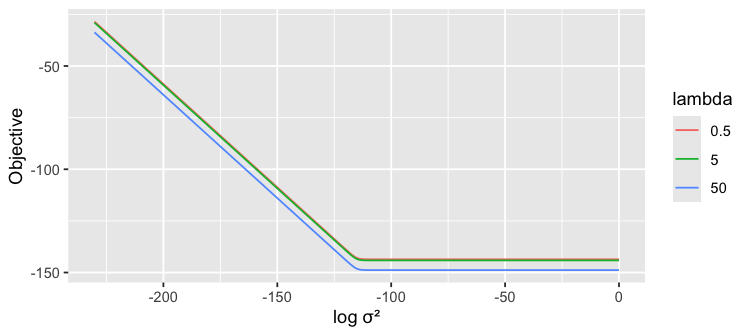}
         \caption{Effect of $\sigma^2$}
         \label{fig:trust2}
     \end{subfigure}
        \caption{Profiles of the objective function $\rho_\lambda^\text{typ}(\mu,\sigma^2)$ for different $\lambda$.}
        \label{fig:trust}
\end{figure}

Lastly, we illustrate the effect of the tuning parameter $\lambda$ in our proposed typicality-based regularization.  Figure~\ref{fig:trust} presents two plots of the objective function based on different $\lambda$.  We postpone investigation into data-driven choices of the tuning parameter for a future paper, but make the observation here that certain values of the tuning parameter make it easier to optimize over $\mu$.  There are, however, still challenges associated with optimizing over $\sigma^2$.  This is to be expected from the fact that the mixture weight $\p$ is very small: at most a few samples out of $n$ will come from the second component of the mixture, hence, the data contains very little information about $\Sigma^2$.  

\subsection{Neyman--Scott problem}
\label{SS:neyman.scott}

Here we take a look at ``a more disturbing example'' \citep{le1990maximum}, due to \citet{neyman.scott.1948}.  In this case, the model $\prob_\theta$, indexed by $\theta=(\xi_1,\ldots,\xi_n,\sigma^2) \in \RR^n \times (0,\infty)$, posits that the data consists of independent pairs, $X_i=(X_{i1},X_{i2})$, where 
\[
X_{i1},X_{i2} \ind \nm(\xi_i, \sigma^2), \quad i=1,...,n.
\]
All of $\Theta=(\Xi_1,\ldots,\Xi_n,\Sigma^2)$ are unknown, but the primary goal is estimation of $\Sigma^2$.  It is easy to show that the maximum likelihood estimator of $\Xi_i$ is $\hat\xi_i = \frac12(X_{i1} + X_{i2})$, and the corresponding maximum likelihood estimator of $\Sigma^2$ is 
\[ \hat\sigma^2 = \frac{1}{2n} \sum_{i=1}^n \{ (X_{i1}-\hat\xi_i)^2 + (X_{i2}-\hat\xi_i)^2 \}. \]
But $\hat\sigma^2$ is a biased estimator---its expected value is $\Sigma^2/2 \neq \Sigma^2$---which further implies that the maximum likelihood estimator is inconsistent.  Hence, the Neyman--Scott problem sheds light on further issues with the principle of maximum likelihood.  

We address this problem by applying the typicality principle, in particular, by employing the typicality-based regularization strategy outlined in Section~\ref{S:first}.  With a slight abuse of the notation in Section~\ref{S:first}, we focus on the interest parameter $\Sigma^2$.  Having eliminated the nuisance parameters $\Xi_1,\ldots,\Xi_n$, a simple goodness-of-fit strategy suggests itself.  That is, we take $\pval_x(\sigma^2)$ to be the (two-tailed) p-value associated with the test statistic 
\[ (x,\sigma^2) \mapsto \frac{1}{\sigma^2} \sum_{i=1}^n \{ (X_{i1}-\hat\xi_i)^2 + (X_{i2} - \hat\xi_i)^2 \} \]
and the corresponding $\chisq(n)$ distribution expected when $\sigma^2$ is the true variance.  Had we chosen not to eliminate the nuisance parameters at the start, we could have employed the same Kolmogorov--Smirnov test procedure discussed above, which amounts to using how closely $\{ (x_{ij} - \xi_i)/\sigma: i=1,\ldots,n; j=1,2\}$ resembles a sample from $\nm(0,1)$ as a measure of how typical $x$ is relative to $(\xi_1,\ldots,\xi_n,\sigma^2)$.  

\begin{figure}[t]
\begin{center}
\includegraphics[width=0.5\linewidth]{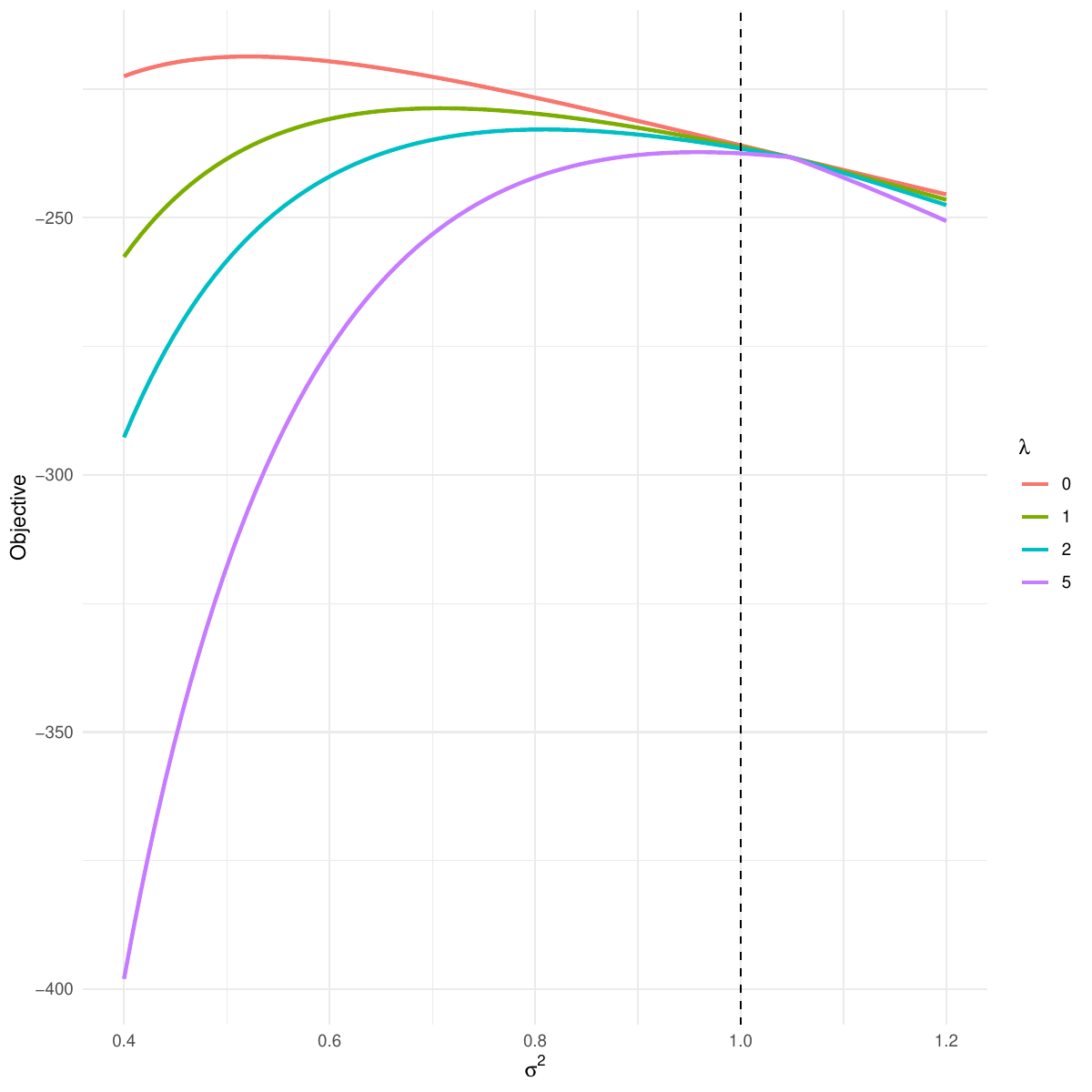}
\caption{
The objective function $\sigma^2 \mapsto \rho_\lambda^\text{typ}(x,\sigma^2)$ with varying choices of $\lambda$. The true $\Sigma^2$ value is denoted as the vertical dashed line.
}
\label{fig:ns}
\end{center}
\end{figure}

For illustration, we simulate data of size $n = 100$, with $\Xi_i \iid N(0,1)$ for $i = 1,\dots,n$ and $\Sigma = 1$. The effect of our proposed regularization is visualized in Figure~\ref{fig:ns}. First note that $\lambda = 0$ corresponds to the usual negative log-likelihood function, which reaches its maximum at approximately $\Sigma^2/2 = 0.5$, consistent with the theoretical analysis. Second, as $\lambda$ increases, the minimizer of the typicality-based objective function gradually shifts from $\Sigma^2/2=0.5$ to $\Sigma^2=1$, suggesting that the proposed estimator is consistent, at least for sufficiently large $\lambda$.  While there is surely room for further investigation, we think it is safe to conclude that our application of the typicality principle resolves the Neyman--Scott paradox and, more general, the shortcomings of the principle of maximum likelihood.

\subsection{Stein's mean vector length}
\label{SS:stein}

Consider the classical problem in which $X$ is a $n$-dimensional normal random vector with unknown mean vector $\Theta$ and identity covariance matrix.  Inference on the mean itself is the same regardless of the dimension, but suppose that the quantity of interest is $\Phi = \|\Theta\|$, the Euclidean length of the mean vector.  Inference on $\Phi$ turns out to be a non-trivial problem, as pointed out by \citet{stein1956, stein1959}, also listed in \citet{fraser.reid.lin.2018} as one of the ``challenge problems'' advanced by the late Sir D.~R.~Cox.  

As a starting point, consider the reparametrization of the model parameter $\Theta$, an $n$-dimensional mean vector, as the pair $(\Phi, \Delta)$, where $\Phi$ is as above and $\Delta = \Theta / \Phi$ is the unit vector direction in which $\Theta$ extends in the $n$-dimensional space.  From this perspective, we see that $\Delta$ is a $(n-1)$-dimensional nuisance parameter so, if $n$ is relatively large, then we can expect challenges with marginalization, similar to those encountered in the Neyman--Scott problem above.  Indeed, Stein demonstrated that the maximum likelihood estimator, $\hat\phi_X = \|X\|$ has non-trivial and non-vanishing upward bias, so the maximum likelihood estimator systematically overestimates $\Phi$.  Two standard strategies for eliminating nuisance parameters involve using log-profile and log-marginal likelihoods which, in this case, are given by 
\[ \ell_x^\text{\sc pro}(\phi) = -\tfrac12( \|x\| - \phi )^2 \quad \text{and} \quad \ell_x^\text{\sc mar}(\phi) = \log q_{\phi^2}(\|x\|^2), \]
where $q_{\phi^2}$ is the density function for a non-central chi-square distribution with $n$ degrees of freedom and non-centrality parameter $\phi^2$.  It is easy to see from the log-profile likelihood function above, and it follows immediately from the familiar invariance property of maximum likelihood estimators, that the maximum profile likelihood estimator of $\Phi$ is the same as the full maximum likelihood estimator, $\hat\phi_X = \|X\|$, which is unsatisfactory as described above.  The log-marginal likelihood maximizer has no closed-form expression, but it closely resembles the simple method of moments estimator $\tilde\phi_X = (\|X\|^2 - n)^{1/2}$, which has a built-in upward bias correction.  

\begin{figure}[t]
\begin{center}
\scalebox{0.7}{\includegraphics{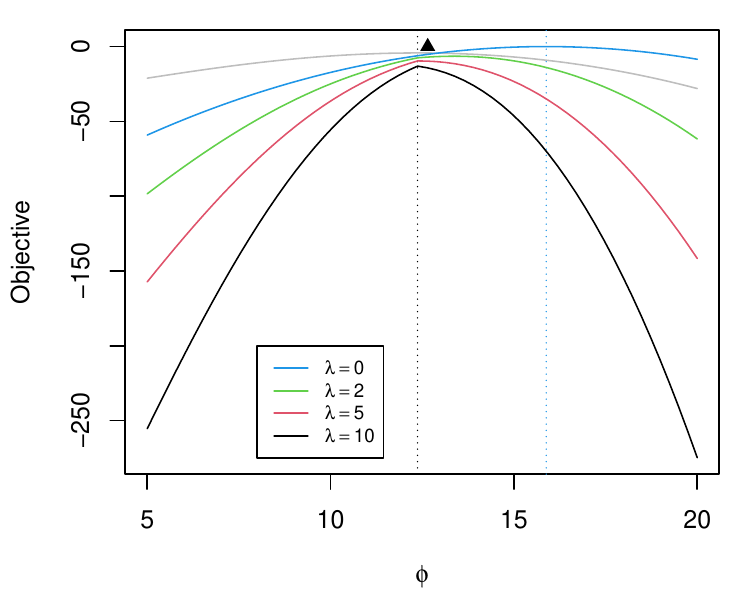}}
\end{center}
\caption{Plots of the objective function \eqref{eq:obj}, specialized for inference on the normal mean vector length in \eqref{eq:obj.stein}, for several different values of the tuning parameter $\lambda \geq 0$.  Gray line is the log marginal likelihood function and the black triangle marks the true value of $\Phi$ in this simulation. The dotted blue and black lines correspond to the objective function maximizers with $\lambda=0$---corresponding to the maximum likelihood estimator---and $\lambda=10$, respectively, the latter clearly being more accurate in this example.}
\label{fig:length}
\end{figure}

As an alternative, after eliminating the nuisance parameters via profiling, we employ our typicality-motivated regularization with objective function 
\begin{equation}
\label{eq:obj.stein}
\rho_\lambda^\text{typ}(x,\phi) =\ell_x^\text{\sc pro}(\phi) -\lambda \, \log \bigl[ \min\{ Q_{\phi^2}(\|x\|^2), 1 - Q_{\phi^2}(\|x\|^2)\} \bigr], 
\end{equation}
where $Q_{\phi^2}$ is the non-central chi-square distribution function corresponding to the density $q_{\phi^2}$ defined above.  Figure~\ref{fig:length} shows plots of this function for several different values of $\lambda$, based on simulated data of size $n$ with true $\Phi = 4 \sqrt{10} \approx 12.65$.  When $\lambda=0$, so that the objective function is the log-profile likelihood, the upward bias is clear, and the maximizer is $\|X\| \approx 15.9$.  As $\lambda$ increases, the peak of the objective function gets sharper and shifts down towards the true $\Phi$; indeed, the objective function maximizer with $\lambda=10$ is $\approx 12.4$.  For comparison, the log-marginal likelihood is shown in gray and, while its maximizer is a satisfactory estimator, it has a rather dull peak indicating a potential loss of efficiency.  

We also carried out a brief simulation study to compare different estimators.  With $\lambda=10$, our estimator has significantly smaller mean squared error than the naive maximum likelihood estimator but very similar mean square error compared to the maximum marginal likelihood estimator.  The Bayesian literature offers a non-informative prior for $\Theta$ designed specifically for inference on $\Phi$ \citep{tibshirani1989}.  This penalizes large $\phi$ values, which makes the corresponding (profiled) maximum {\em a posteriori} estimator appealing.  But this prior has a strong built-in ``$\lambda$,'' which is $O(n)$.  So it happens that this approach penalizes too much, pushing the estimates too close to 0, hence is not competitive with our estimator or the usual marginal maximum likelihood estimator.

\section{Reliable uncertainty quantification}
\label{S:broader}

\subsection{Typicality principle, revisited}
\label{SS:second}

Of course, there's more to statistics and data science than point estimation, and it turns out that the typicality principle has implications beyond the important but relatively narrow point estimation context in Section~\ref{S:first}.  We start with a discussion of the importance and novelty of our direct and specific focus on typicality compared to other approaches that have been taken previously.  

When prior information about $\Theta$ is vacuous, it is mathematically impossible to derive a posterior distribution for $\Theta$, given $x$, in a rigorous, justifiable way.  With only data $x$ and model $\{\prob_\theta: \theta \in \TT\}$, all we can possibly say about a hypothesis $H$ is if $x$ is typical or atypical relative to $\prob_\theta$ for $\theta \in H$.  Therefore, any legitimate attempt to quantify uncertainty about $\Theta$ given $x$ is simply choosing to interpret typicality---a characteristic of data given a value of the parameter---as a statement about that value of the parameter given data.  This is what's at the heart of falsificationism and, specifically, of Mayo's
\begin{quote}
Frequentist Principle of Evidence: {\em Drawing inferences from data requires considering the relevant error probabilities associated with the underlying data generation process} \citep{mayo2014}.
\end{quote}
Neyman's frequentism stops with the intuitive version of the typicality principle: 
\begin{itemize}
\item specify a test statistic, a significance level, and a rejection region such that, if hypothesis $H$ is true, then the event ``test statistic falls in the rejection region'' has (error) probability no more than the stated significance level; and
\vspace{-2mm}
\item say data $x$ is sufficiently atypical relative to $H$ if the test statistic based on $x$ falls in the stated rejection region.
\end{itemize} 
That is, pure frequentists make no attempt to quantify uncertainty about the truthfulness of $H$ given $x$---they're satisfied with an error-probability-controlling decision rule.  Attempts to find the Bayesian--frequentist Holy Grail go further by constructing data-dependent measures of support for or degrees of belief/confidence in the truthfulness of hypothesis $H$ that remain consistent with Mayo's frequentist principle of evidence.  This list of attempts includes default-prior Bayesian inference \citep[e.g.,][]{jeffreys1946, bergerbernardosun2009, datta.ghosh.1995}, Fisher's fiducial inference \citep[e.g.,][]{fisher1933, fisher1935a, zabell1992} and its generalizations \citep[e.g.,][]{fraser1968, hannig.review, xie.singh.2012}, Dempster--Shafer theory \citep[e.g.,][]{dempster1966, dempster2008, shafer1976, shafer1982}, and inferential models \citep[e.g.,][]{imbasics, imbook, plausfn, imchar, martin.basu}.  One thing these approaches have in common is that, at least at face value, they fail to recognize that typicality as the sole fundamental notion to work with; this failure creates opportunities for confusion and prevents progress.  By putting typicality front and center, we are able to clear up this confusion.  In particular, we show what kinds of properties are consistent with a measure of typicality and, in turn, demonstrate what it takes to find the Holy Grail of reliable and principled uncertainty quantification about unknowns.  

In order to help orient our pursuit, we start by giving a formal counterpart to the informally stated typicality principle in Section~\ref{S:first}.  We follow this up with a clarification of several terms and concepts introduced in the formal statement.  

\begin{tprin}[formal version]
Given a level $\alpha \in (0,1)$ and typicality measure $\tau_x: 2^\TT \to [0,1]$ that is large when $x$ is consistent with $H$ and small otherwise, if data $x$ is sufficiently atypical relative to the posited model and hypothesis $H$, in the sense that $\tau_x(H) \leq \alpha$, then hypothesis $H$ is unwarranted at level $\alpha$. 
\end{tprin}

There are two key features inherent in a measure of typicality that will assist in making the abstract principle above more concrete.  First, typicality of $x$ relative to a general, non-singleton $H$ has no direct meaning---the typicality of $x$ relative to $H$ can only be assessed indirectly via the typicality of $x$ relative to $\theta$ for individual $\theta$ in $H$.  That is, one must first define typicality of $x$ relative to simple, singleton hypotheses $H=\{\theta\}$, for $\theta \in \TT$, and then derive the typicality $\tau_x(H)$ of $x$ relative to general $H$ from the primitive typicality $\tau_x(\theta) = \tau_x(\{\theta\})$.  An important observation is that, unlike a Bayesian posterior for a continuous parameter that assigns zero posterior probability to ``$\Theta=\theta$'' for all $\theta$, the typicality $\tau_x(\theta)$ of $x$ relative to $\theta$ is generally non-zero.  To remain consistent with our falsificationist motivations, it is essential that a hypothesis be falsified if and only if all those hypotheses that entail it are falsified.  This means data is sufficiently atypical relative to ``$\Theta \in H$'' if and only if it's sufficiently atypical relative to ``$\Theta=\theta$'' for all $\theta \in H$.  Mathematically, it means that $\tau_x(H)$ cannot exceed $\sup_{\theta \in H} \tau_x(\theta)$ since, otherwise, it would be possible for $x$ to be typical relative to $H$ but sufficiently atypical relative to all $\theta$ in $H$.  Since $\tau_x(H)$ for non-singleton hypotheses $H$ is derived from the primitive $\tau_x(\theta)$, and since there is no practical advantage to taking $\tau_x(H)$ smaller than necessary, we take 
\begin{equation}
\label{eq:maxitive}
\tau_x(H) = \sup_{\theta \in H} \tau_x(\theta), \quad H \subseteq \TT. 
\end{equation} 
In other words, $H \mapsto \tau_x(H)$ is a maxitive set function \citep[][Ch.~1.2]{molchanov2005}, just like the plausibility function associated with a consonant belief function \citep[][Ch.~10]{shafer1976} or a possibility measure \citep{dubois.prade.book}.  

Second, since what's being measured is typicality, a suitable calibration property is expected.  If $\Theta=\theta$, then we expect the random variable $\tau_X(\theta)$ to be moderate to large but, of course, it's possible that $\tau_X(\theta)$ is small, and this latter event should likewise have small probability.  So, when we say that ``$x$ is atypical relative to $\theta$,'' so that $\tau_x(\theta)$ is small, we mean that a $\prob_\theta$-rare has occurred.  Then it makes sense to require that $\tau_X(\theta)$ have a $\unif(0,1)$ distribution, or at least be stochastically no smaller than $\unif(0,1)$, as a function of $X \sim \prob_\theta$, i.e., 
\begin{equation}
\label{eq:valid}
\sup_{\theta \in \TT} \prob_\theta\{ \tau_X(\theta) \leq \alpha \} \leq \alpha, \quad \alpha \in [0,1]. 
\end{equation} 
The above property could be generalized to $\sup_\theta \{ \tau_X(\theta) \leq g(\alpha) \} \leq \alpha$ for some known, strictly increasing function $g: [0,1] \to [0,1]$.  It's necessary that the function $g$ be known because, otherwise, the scale of our typicality measure would be likewise unknown and we would be unable to label $x$ as typical/atypical.  Importantly, this known function $g$ also cannot depend on $\theta$, $X$, or any other potentially relevant features of the data-generating process, such as sample size.  To explain this we make an analogy to the interpretation of weather forecasting: a claim ``10\% chance of precipitation today'' means the same regardless of the kind of precipitation in question, the kind of climate model is used to generate the prediction, or the kind/amount of data used.  So, if the function $g$ above must be known and independent of $\theta$ and $X$, then the desired property \eqref{eq:valid} is without loss of generality since, if necessary, we could just redefine $\tau_x(\theta)$ as $g^{-1} \tau_x(\theta)$ and obtain an equivalent typicality measure for which \eqref{eq:valid} holds.  

Putting everything together, under the typicality principle with a typicality measure satisfying the above properties, there are a number of desirable inference-related consequences.  Two specific consequences follow immediately from the discussion above.
 
\begin{prop}
A test procedure that rejects the hypothesis $H$ if and only if $x$ is sufficiently atypical relative to $H$ at level $\alpha$, i.e., rejects if and only if $\tau_x(H) \leq \alpha$, controls the frequentist Type~I error in the sense that 
\begin{equation}
\label{eq:valid.H}
\sup_{\theta \in H} \prob_\theta\{ \tau_X(H) \leq \alpha \} \leq \alpha, \quad \alpha \in [0,1], \quad H \subseteq \TT. 
\end{equation}
\end{prop}

\begin{prop}
\label{prop:conf.set}
The set $C_\alpha(x) = \{\theta \in \TT: \tau_x(\theta) \geq \alpha\}$, i.e., the collection of all $\theta$'s such that $x$ is not sufficiently atypical relative to $\theta$ at level $\alpha$, is a nominal $100(1-\alpha)$\% confidence set for $\Theta$ in the sense that 
\[ \sup_{\theta \in \TT} \prob_\theta\{ C_\alpha(X) \not\ni \theta \} \leq \alpha, \quad \alpha \in [0,1]. \]
\end{prop} 

These are exactly the kind of frequentist error rate control guarantees that ensure basic reliability of inference.  Connections to point estimation will be drawn in Section~\ref{SS:practice} once we discuss construction of typicality measures. The property in \eqref{eq:valid.H} holds for all fixed $H \subseteq \TT$, but reliable uncertainty quantification requires that $H$ can vary and that similar reliability properties are maintained.  The details are a bit beyond our present scope but it is shown in \citet{cella.martin.probing} that a uniform-in-$H$ version of \eqref{eq:valid.H} holds, which means the typicality measure can be used for reliable probing of various hypotheses.  Finally, since $\tau_x$ has at least some loose resemblance to a Bayesian or fiducial posterior distribution \citep[cf.,][]{martin.isipta2023}, it is natural to consider formal decision-making wherein the ``optimal'' action is defined as that which minimizes a suitable upper expected loss---a Choquet integral \citep[e.g.,][App.~C]{lower.previsions.book}---with respect to $\tau_x$.  \citet{imdec} showed that the frequentist calibration properties inherent in $\tau_x$ carry over in novel ways to the IM's assessment of the $\Theta$-dependent, to-be-incurred loss associated with an action.  Again, the details are beyond our present scope.

\subsection{Putting the principle into practice}
\label{SS:practice}

In light of the developments discussed in Section~\ref{S:first}, there is a conceptually straightforward path to getting from the typicality principle---both the intuitive and formal versions---to a practical typicality measure that has the properties discussed in the previous section.  Let $\check\theta_x = \check\theta(x)$ denote the minimizer of the objective function $\theta \mapsto \rho_\lambda^\text{typ}(x,\theta)$ defined in \eqref{eq:obj}, which depends implicitly on the choice of $\lambda \geq 0$, with $\lambda=0$ leading to the usual maximum likelihood estimator.  Now define 
\begin{equation}
\label{eq:tau}
\tau_x(\theta) = \prob_\theta \bigl\{ R_\lambda(X,\theta) \geq R_\lambda(x, \theta) \bigr\}, \quad \theta \in \TT, 
\end{equation} 
where 
\begin{equation}
\label{eq:rl}
R_\lambda(x,\theta) = \rho_\lambda^\text{typ}(x,\theta) - \rho_\lambda^\text{typ}(x,\check\theta_x). 
\end{equation}
While the ingredients might be unfamiliar, the particular operation is a familiar one: it's the p-value associated with a test of the hypothesis ``$\Theta=\theta$'' using data $x$ and test statistic $R_\lambda(x,\theta)$.  This can also be seen as a version of the probability-to-possibility transform \citep[e.g.,][]{hose2022thesis, dubois.etal.2004}, critical to the statistical developments in \citet{martin.partial2}.  Below we briefly discuss the properties of the $\tau_x$ in \eqref{eq:tau}.  

Since the right-hand side of \eqref{eq:tau} is a probability, the left-hand side clear cannot exceed 1.  But it actually attains the value 1 at $\theta=\check\theta_x$, the estimator from Section~\ref{S:first}.  Therefore, $\check\theta_x$ can legitimately be called the ``maximum typicality estimator,'' and is further distinguished as being the point contained in all of the confidence sets $C_\alpha(x)$ as $\alpha$ ranges from 0 to 1.  That $\theta \mapsto \tau_x(\theta)$ attains the value 1 further implies that the typicality measure $H \mapsto \tau_x(H)$ as in \eqref{eq:maxitive} is a genuine possibility measure and, therefore, can be interpreted as a coherent upper probability.  This is in addition to all the desirable frequentist properties in Section~\ref{SS:second} above.  Is this the aforementioned Bayesian--frequentist Holy Grail?

Implementation of the proposed framework requires that we can evaluate the contour defined in \eqref{eq:tau}.  There's a simple and easy-to-explain approach that we describe next, which is sufficient for most problems involving low-dimensional parameters.  Fix an argument $\theta$ to the typicality contour, and then approximate \eqref{eq:tau} as 
\begin{equation}
\label{eq:tau.naive}
\tau_x(\theta) \approx \frac1M \sum_{m=1}^M 1\{ R_\lambda(X_{m,\theta}, \theta) \geq R_\lambda(x, \theta) \}, \quad \theta \in \TT,
\end{equation}
where $X_{m,\theta}$ are independent copies of the data $X$ drawn from $\prob_\theta$, for $m=1,\ldots,M$.  It may be the case that the above approach is too expensive for practical use.  Indeed, in moderate- to high-dimensional problems, having to evaluate the contour on a sufficiently dense grid covering the relevant portion of the parameter space requires substantial computational investment.  Various adjustments can be made, e.g., using importance sampling, to reduce the burden of generating such a diverse set of Monte Carlo samples, but the extent of these improvements is limited.  For this reason, recent efforts have focused on the development of new strategies that mimic the Monte Carlo methods used by Bayesians, i.e., where samples of the parameter values---not new data sets---are drawn from a ``posterior distribution'' where the curse of dimensionality can be better controlled.  These details are beyond the scope of the present paper, but the interested reader can consult \citet{calibrated.boostrap} and \citet{immc}.

\subsection{Relation to other statistical principles}
\label{SS:relation}


The most familiar statistical principle is the {\em likelihood principle} \citep[e.g.,][]{birnbaum1962, basu1975, bergerwolpert1984}, which states that all of what is relevant in the data for inference on $\Theta$ is captured by the shape of the likelihood function.  This might seem intuitive and harmless, since the commonly used maximum likelihood estimators and likelihood ratio statistics only depend on the shape of the likelihood function.  On closer inspection, however, the things we commonly do with these summaries, e.g., p-value calculations, rely on sampling distributions under the posited model and, since sampling distributions aren't determined by the observed likelihood function, inference based on these violates the likelihood principle.  On its own, this violation is of no concern, but becomes potentially problematic in light of Birnbaum's theorem stating that the likelihood principle is equivalent to conjunction of the more commonsense sufficiency and conditionality principles; accordingly, a violation of the likelihood principle implies a violation of at least one commonsense principle, hence the controversy.  But legitimate doubts about the reach of Birnbaum's theorem, first back in \citet{durbin1970} and more recently in \citet{evans2013} and \citet{mayo2014}, have only fueled the controversy.  

We emphasize 
that the typicality principle is unapologetically orthogonal to the likelihood principle.  It is clear, for example, that unless $\lambda=0$, the minimizer $\check\theta(x)$ of the objective function \eqref{eq:obj} is not determined by the shape of the likelihood alone.  At a fundamental level, the point is that data being ``(a)typical'' can only be assessed relative to the posited sampling model, which, again, is not determined by the observed likelihood function.  If one expects their inferences to favor models that adequately predict the observed data, then observed likelihood alone simply isn't enough.  

From the perspective of Sections~\ref{SS:second}--\ref{SS:practice} above, there is more that can be said even if we ignore the typicality-motivated penality term $r_\lambda^\text{typ}$ or, equivalently, if we take $\lambda=0$ in \eqref{eq:obj}.  In that case, $R_0(x,\theta)$ in \eqref{eq:rl} is just the log-relative likelihood which, itself, only depends on the shape of the observed likelihood function.  It is easy to verify that inference based on the typicality measure \eqref{eq:tau} satisfies the sufficiency principle.  It can also readily facilitate conditioning on an ancillary statistic, if the user so desires, but we refer the reader to \citet[][Sec.~6]{martin.partial2} for these details.  But despite the close connection to the relative likelihood, inference based on the typicality measure \eqref{eq:tau} still violates the likelihood principle.  The outer $\prob_\theta$-probability calculation is the culprit, as this brings in aspects of the model not included in the relative likelihood itself.  As \citet{martin.basu} argues, however, the modified version of the typicality measure \eqref{eq:tau}, 
\begin{equation}
\label{eq:tau.lp}
\tau_x^\text{\sc lp}(\theta) := \sup \prob_\theta\{ R_0(X,\theta) \geq R_0(x,\theta) \}, \quad \theta \in \TT, 
\end{equation}
where the outer supremum is over all those probabilities $\prob_\theta$ that admit a relative likelihood with the same shape as $R_0(x,\theta)$, retains the above frequentist properties and simultaneously satisfies the likelihood principle.  (Practically, the above supremum is over different stopping rules, procedures for deciding when to stop collecting data and begin the analysis; the effect of the stopping rule cancels in the likelihood ratio, but clearly each stopping rule determines its own unique probability model.)  The above modification makes crystal clear the cost of satisfying the likelihood principle: since $\tau_x(\theta) \leq \tau_x^\text{\sc lp}(\theta)$ for all $x$ and $\theta$, usually with strict inequality ``$\ll$,'' and since smaller typicality measure means more efficient inference, the cost of satisfying the likelihood principle is a potentially significant loss of efficiency.  If we genuinely don't know the sampling model (or stopping rule), then the modification in \eqref{eq:tau.lp} is justified.  If we are confident in the posited sampling model (or stopping rule), then there is no need to sacrifice efficiency in order to satisfy the likelihood principle; this is why we're ``unapologetic'' about our proposal violating the likelihood principle.  There is a middle-ground between the two extremes where, e.g., we know that one in a proper subset of all possible stopping rules was used, and we refer the reader to \citet{martin.basu} for details on how this can be accommodated.  

Another statistical principle, relevant to our proposal here, was mentioned in Section~\ref{S:intro}, namely, the {\em prediction principle} advanced in \citet{sts.discuss.2014, imbook}.  This principle states that, roughly, the sampling model should be expressed in terms of a ``predictable quantity'' and that inference should be based on (reliable) prediction of that predictable quantity.  It turns out that this is closely related to the proposed typicality principle, as we describe next.  The unobserved realization of a random variable $U$ with known distribution $\prior$ supported on $\UU$ is an example of a predictable quantity.  Martin and Liu suggest predicting the unobserved realization of $U$ as follows: define a set-valued function $\S: \UU \to 2^\UU$ such that the corresponding hitting probability function 
\[ h_\S(u) = \prior\{ \S(U) \ni u \}, \quad u \in \UU, \]
has the property 
\begin{equation}
\label{eq:valid.U}
\prior\{ h_\S(U) \leq \alpha \} \leq \alpha, \quad \alpha \in [0,1]. 
\end{equation}
A relevant example is where $\UU=[0,1]$, $\prior=\unif(0,1)$, and $\S(u) = [0,u]$, so that $h_\S(u) = 1-u$ and the property \eqref{eq:valid.U} holds trivially.  If the goal is to guess an unobserved realization of a random variable drawn from $\prior$, then we have no hope but to assume that it's a ``typical'' realization, and we measure the typicality of a candidate $u$ via $h_\S(u)$.  The intuition is that, by \eqref{eq:valid.U}, the random set $\S(U)$, with $U \sim \prior$, contains typical realizations from $\prior$ in the sense that the random hitting probability $h_\S(U)$ tends to be not-small; hence we're warranted to conclude that $u$ values with $h_\S(u)$ small are atypical and, hence, poor predictions of the unobserved realization.  

The connection to our discussion of typicality in the context of statistical inference is as follows.  \citet{plausfn, gim} suggests a generalized association to relate observed data $x$, uncertain parameter $\theta$, and predictable quantity $u$.  Based on the formulation in Section~\ref{SS:practice}, a natural choice of this generalized association is 
\[ G_\theta\{ R_\lambda(x,\theta) \} = u, \quad (x,\theta,u) \in \XX \times \TT \times \UU, \]
where $G_\theta$ is the distribution function of the random variable $R_\lambda(X,\theta)$, as a function of $X \sim \prob_\theta$, and $\UU=[0,1]$. Using the same set-valued map in the illustration above, the typicality measure in \eqref{eq:tau} can be re-expressed as 
\[ \tau_x(\theta) = h_\S\bigl[ G_\theta\{ R_\lambda(x,\theta) \} \bigr], \quad \theta \in \TT, \]
so that our proposed notion of typicality here can be directly tied to the prediction of predictable quantities emphasized in \citet{sts.discuss.2014}.

\subsection{Stein's mean vector length, again}
\label{S:more-examples}

To illustrate the broader, typicality-motivated uncertainty quantification strategy described in Section~\ref{S:broader}, we revisit Stein's normal mean vector length example from Section~\ref{SS:stein} above.  Again, suppose we have data $X=(X_1,\ldots,X_n)$ consisting of independent Gaussian random variables where $X_i \sim \nm(\Theta_i, 1)$ with $\Theta_i$ uncertain, for $i=1,\ldots,n$, and the goal is to quantify uncertainty about the length $\Phi = \|\Theta\|$ of the mean vector $\Theta=(\Theta_1,\ldots,\Theta_n)$.  More specifically, with a slight abuse of notation and terminology, our goal is to construct the typicality contour $\phi \mapsto \tau_x(\phi)$ for the length $\Phi$ of the mean vector $\Theta$.  Towards this, a key observation is that the proposed objective function $\phi \mapsto \rho_\lambda^\text{typ}(X,\phi)$ and, hence, the corresponding $R_\lambda(X,\phi)$ in \eqref{eq:rl} depends on data $X$ only through the squared length $\|X\|^2$.  Consequently, the probability calculation in the general formula \eqref{eq:tau} here only depends on $\phi$, since the distribution of $\|X\|^2$ under the posited model is the non-central chi-square with non-centrality parameter $\phi^2$, so dependence on the nuisance parameter---the mean vector direction---drops out automatically.  While elimination of the nuisance parameter is relatively straightforward in this case, it will not be so straightforward in other cases; see, e.g., \citet{basu1977}, \citet{immarg}, and \citet{martin.partial3} for further discussion of marginalization.  

\begin{figure}[t]
\begin{center}
\scalebox{0.7}{\includegraphics{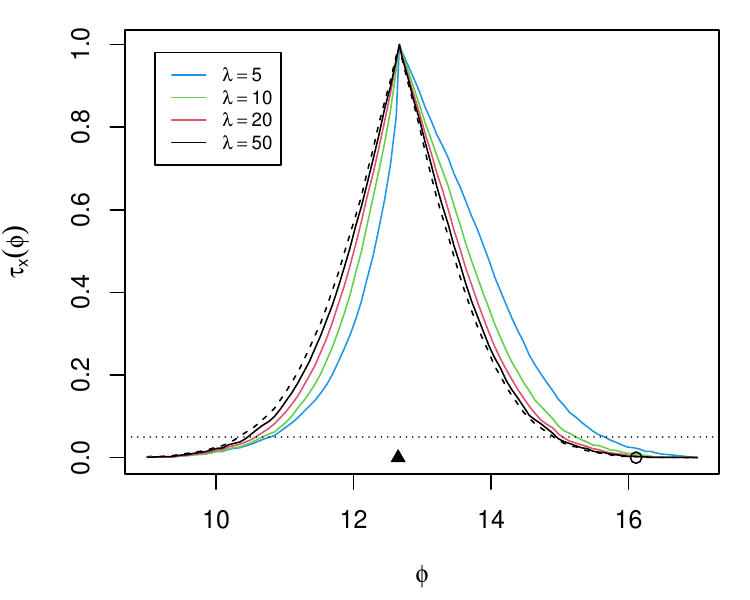}}
\end{center}
\caption{Plots of the typicality contour $\phi \mapsto \tau_x(\phi)$ for the mean vector length $\Phi$ based on the simulated data of size $n=100$ and different $\lambda$ values. Dashed line corresponds to the marginal likelihood-based contour suggested in \citet{martin.partial3}. Open circle corresponds to the maximum likelihood estimator $\|x\|$ and the triangle is the true $\Phi$.}
\label{fig:stein.pl}
\end{figure}

For illustration, we simulated data with $n=100$ and $\Theta$ such that $\Phi=4\sqrt{10} \approx 12.65$; the generated data $X$ is such that $\|X\| \approx 16.10$.  Figure~\ref{fig:stein.pl} shows a plot of the typicality contour $\tau_x(\phi)$ as defined in \eqref{eq:tau}, evaluated using the Monte Carlo approximation \eqref{eq:tau.naive} described in Section~\ref{SS:practice} above for several different values of $\lambda$.  The first observation is that the contours indexed by $\lambda$ are all peaked at effectively the same point, which is the typicality-motivated estimator of $\Phi$ proposed in Section~\ref{S:first}.  As is common of these contour functions, they vanish away from the peak, often asymmetrically.  The dotted horizontal line at 0.05 determines the 95\% confidence interval for $\Phi$ as in Proposition~\ref{prop:conf.set}; note that the maximum likelihood estimator does not belong to these confidence sets in this case, but the true $\Phi$ does.  For comparison, we have also shown the contour corresponding to the recommended marginal likelihood-based construction in \citet{martin.partial3}; that is, the $R_\lambda(x,\phi)$ as in \eqref{eq:rl} is taken to be the negative log-marginal likelihood, based on the non-central chi-square distribution of $\|X\|^2$.  While there is no agreed-upon ``best'' solution to this marginal inference problem, at least to our knowledge, it is difficult to imagine how another solution can be substantially better than that based on the marginal likelihood: some information about $\Phi$ is lost when the data $X$ is reduced to $\|X\|^2$, but profiling or integrating with respect to a prior distribution will introduce bias that needs to be corrected for.  An interesting observation, although not unexpected given the penalty structure in \eqref{eq:obj.stein}, is that the proposed typicality contour merges with the high-quality marginal likelihood-based solution in \citet{martin.partial3} as $\lambda$ increases.  


\section{Conclusion}
\label{S:discuss}

Motivated by the deep philosophical considerations and scientific attitudes of Popper and Tukey, here we advanced a new {\em typicality principle} that has a number of methodological and foundational implications for statistics and data science.  First, on the methodological side, the typicality principle immediately suggests a novel strategy for regularization in the context of parameter estimation.  Indeed, instead of shrinking estimators toward assumed structure (e.g., ``sparsity'') in the estimand, our typicality-based focus encourages goodness-of-fit, ensuring that the observed data looks ``typical'' under the fitted model.  Three non-trivial, illustrative examples are presented where the method of maximum likelihood fails miserably while our proposed typicality-focused regularization strategy is shown to be more than satisfactory.  Second, on the foundational side, a more formal version of the typicality principle is readily accommodated by the general inferential model framework for provably reliable uncertainty quantification beyond point estimation, hypothesis testing, etc.  This broader form of uncertainty quantification can readily accommodate von Neumann and Morganstern-style decision-making and other formal inferences, similar to Bayes, but without requiring a prior distribution and without sacrificing on error-rate control guarantees.  More generally, we believe that various implementations of the proposed typicality principle will be beneficial to data science, as automated applications are expected to play a key role in the advancement of artificial intelligence, etc.  Although our focus in the present paper was on model-based inference, the very notion of typicality is tied tightly to goodness-of-fit considerations, so we fully expect the typicality principle and the methodology derived from it to have an impact on scientific modeling as well as inference.  

The typicality principle and the various methodological advancements derived from it are subject to further investigations, applied and theoretical.  On the application side, modern problems in data science involve complex models where, without the guidance of one regularization strategy or another, the tendency will be to overfit; it is precisely this tendency to overfit that leads to the poor performance of maximum likelihood estimators in the examples presented in Section~\ref{S:applications}.  Then a natural follow-up to the present paper is an investigation into the performance of the proposed typicality-based regularization strategy in a class of modern, data science-related problems involving complex and over-parameterized models such as deepnets and transformers \citep{Vaswani2017Attention}. It would also be interesting to compare our proposed with other advanced techniques, such as knowledge distillation \citep[see, e.g.,][]{hinton2015,jiang2022estimation}. There is also the important practical question of how to set the tuning parameter $\lambda$ in \eqref{eq:obj}.  While there are so many now-standard tuning parameter selection strategies available, a relevant question is if the data-dependence inherent in our typicality-based penalty warrants new tuning parameter selection considerations.  After all, compared to the usual sparsity-encouraging penalties, p-values have a meaningful scale, so new considerations concerning how to balance the contribution of the latter penalty with the likelihood may be needed.  On the theoretical side, the finite- and large-sample efficiency properties of the proposed typicality-based regularized maximum likelihood estimators, and of the broader uncertainty quantification developed in Section~\ref{S:broader}, are completely open for investigation.  Nevertheless, when necessary, the familiar sparsity-based penalties can be easily incorporated to offer a hybrid regularization approach.

Aligned with our proposed statistical principles, and motivated by other deep philosophical considerations, recent advancements have been made to enhance the creativity and trustworthiness of artificial intelligence \citep[e.g.,][]{eschker2024}. An major challenge lies in understanding how these and other philosophical advances can help to refine today's cutting-edge methods and inspire new developments that push modern boundaries.  It is also important that these philosophical contributions be of Tukey's hands-dirty, applications-oriented, ``bottom-up'' style, as opposed to hands-tying, ``top-down'' style protocol dictated in an ivory tower.

\section*{Acknowledgement}

Liu and Zhang are supported by the U.S. National Science Foundation grant DMS-2412629. Martin is supported by the U.S. National Science Foundation grants SES–2051225 and DMS–2412628.


\bibliographystyle{apalike}
\bibliography{mybib}

\end{document}